\documentclass{amsart}
\usepackage{amsthm}
\usepackage{amsfonts}
\usepackage{amsmath}
\usepackage{graphicx}
\usepackage{float} 
\newtheorem{thm}{Theorem}[section]

\newtheorem{coro}[thm]{Corollary}

\begin{document}

\title{ construction of Steiner Triple Systems of type $v\longrightarrow 2v+7$}

\author{Paola Bonacini}
\address{Dipartimento di Matematica e Informatica\\
Viale A Doria 6,
95125 Catania, Italy}
\email{paola.bonacini@unict.it}

\author{Mario Gionfriddo} 
\address{Dipartimento di Matematica e Informatica\\
Viale A Doria 6,
95125 Catania, Italy}
\email{mario.gionfriddo@unict.it}

\author{Lucia Marino}
\address{Dipartimento di Matematica e Informatica\\
Viale A Doria 6,
95125 Catania, Italy}
\email{lucia.marino@unict.it}

\begin{abstract}
    
A \textit{Steiner Triple System} ($STS$) of order $v$ is a hypergraph uniform of rank 3, with $v$ vertices and such that every 2-subset of vertices has degree 1. In this paper we give a construction, by difference method, of type $v\longrightarrow 2v+7$ with $v=2^n-7$, which means that, given an $STS$ of order $v=2^n -7$, it is always possible to construct an $STS$ of order $2^{n+1}-7$. Through this construction it is possible to get for any $n\ge 5$ an $STS(2^n-7)$ with a maximal independent set of maximal cardinality and which is $(n-1)$-bicolorable.
\end{abstract}


\maketitle

\section{INTRODUCTION}
A Steiner system $S(h,k,v)$ is a pair $\Sigma=(X,\mathcal{B})$, where $X$ is a $v$-set and $\mathcal{B}$ is a family of $k$-subsets of $X$ such that every $h$-subset of $X$ is contained in exactly one member of $\mathcal{B}$ (see \cite{CR,GMV,LR,LR2}). Using hypergraph theory terminology, a Steiner system is a hypergraph $\Sigma=(X,\mathcal{B})$ of order $v$, uniform of rank $k$, such that every $h$-subset $Y$ of $X$ has degree $d(Y)=1$(see \cite{GMV}).

A Steiner Triple System ($STS$) is a system $S(2,3,v)$. Steiner systems $S(h,k,v)$ were defined for the first time by Woolhouse in 1844 \cite{W}, who asked for which positive integers $h,k,v$ an $S(h,k,v)$ there exists. This problem remains unsolved in general until today, even if many partial results have been given. In 1847 T. Kirkman \cite{K} and J. Steiner \cite{S}, independently, showed that \textit{an STS(v) there exists if and only if} $v\equiv 1\hspace{0.5ex}$or$\hspace{0.5ex} 3\hspace{0.5ex} \mod 6$.

Other results have been determined by H. Hanani about the spectrum of $S(3,4,v)$ and $S(2,4,v)$, respectively in 1960 \cite{H} and in 1962 \cite{H2}. In the literature there are many constructions to obtain an $STS$, starting from a given $STS(v)$. Among them, it is the well known construction, indicated by $v\longrightarrow 2v+1$, which gives an $STS(2v+1)$ starting from an $STS(v)$. Other constructions are of type $v\longrightarrow 3v$.

Given an $STS$ $\Sigma=(X,\mathcal B)$, an \emph{independent set} $T$ is a subset $T\subset X$ that doesn't contain any triple of $\mathcal B$. $T$ is a \emph{maximal independent set} for $\Sigma$ if it is not contained in another independent set. By \cite{SS} a maximal independent set in an $STS(v)$, with $v\equiv 1,9\mod 12$, has cardinality at most $\tfrac{v-1}{2}$ (see also \cite{CD}). Note that $v=2^n-7\equiv 1\text{ or }9\mod 12$ for any $n\in\mathbb N$.

A \emph{bicoloring} of an STS $\Sigma=(X,\mathcal B)$ is a coloring of the vertices in such a way that all the blocks contain vertices of exactly two colors. $\Sigma$ is called $k$-bicolorable if there exists a bicoloring of $\Sigma$ using exactly $k$ colors. For a bicolorable Steiner Triple System $\Sigma$ we denote $\chi(\Sigma)$ (resp. $\overline{\chi}(\Sigma)$) the \emph{lower} (resp. \emph{upper}) chromatic number, which is the smallest (resp. largest) integer $k$ for which there exists a $k$-bicoloring of $\Sigma$ using exactly $k$ colors. Moreover, by a simple count it is immediate to see that no $STS$ is $2$-bicolorable (see \cite{R}).

These type of colorings were introduced by Voloshin \cite{V1,V2} in the context of mixed hypergraphs and for Steiner systems the study was initiated by Milazzo and Tuza \cite{MT}. It is relevant to remark that in \cite{MT} it was proved that there exists a strong connection between bicolorable $STS(2^n-1)$ with the highest possible upper chromatic number and the ``doubling plus one construction'' $v\mapsto 2v+1$. Indeed, if $\Sigma$ is an $STS(2^n-1)$ with $\overline{\chi}(\Sigma)=n$, then $\Sigma$ is obtained from the $STS(3)$ by repeated applications of  ``doubling plus one constructions''. Later, many other papers dealt with bicolorings of Steiner systems (see, for example, \cite{BGGMTV,BGMV,CDR,GGM,MTV}).

In this paper we give in Theorem \ref{T} a construction $v\longrightarrow 2v+7$, in which $v=2^n-7$. We show that, through an iteration of this construction, it is possible, starting from the $STS(9)$, to obtain for any $n\ge 5$ an $STS(2^n-7)$, $\Sigma$, with a maximal independent set of maximal cardinality and with $\overline{\chi}(\Sigma)=n-1$, where $n-1$ is the largest possible upper chromatic number.

\section{The well-known construction $v\longrightarrow 2v+1$}

It is well-known that it is always possible to construct an $STS(2v+1)$ starting from an $STS(v)$.

\begin{thm} 
    If $\Sigma=(X,\mathcal{B})$ is an $STS(v)$, then there exists an $STS(2v+1)$ embedding $\Sigma$.
\end{thm}
\begin{proof}
Let $\Sigma=(X,\mathcal{B})$ be an STS(v) defined on $X=\{x_1, x_2, .... , x_v\}$ Further, let:
\begin{itemize}
    \item $Y=\{y_1, y_2, .... , y_{v+1}\}$ be a set of cardinality $v+1$ (even number) such that $X\cap Y=\emptyset$;
    \item $\mathcal{F}=\{F_1, F_2, .... , F_{v}\}$ be a factorization of the complete graph $K_{v+1}$ defined on $Y$;
\item $\varphi: X\longrightarrow \mathcal{F}$ be any bijection from $X$ into $\mathcal{F}$ .
\end{itemize}
Define the hypergraph $\Sigma'=(X',\mathcal{B}')$ as follows:
\begin{itemize}
    \item $X'=X\cup Y$;
    \item $\mathcal{B}'=\mathcal{B}\cup \Gamma$, where:
$$ \Gamma=\{ \{x,y,z\}: x\in X, {y,z}\in \varphi(x)\}.$$
\end{itemize}
We call the triples of $\mathcal{B}$ \textit{of type 1}, the triples of $\Gamma$ \textit{of type 2}. We say that $\Sigma'$ is an $STS(2v+1)$. Indeed:
\begin{enumerate}
    \item  it is immediate that $\Sigma'$ is a hypergraph of order $|X'|=2v+1$, uniform of rank $3$;
    \item for every $x,y\in X', x\neq y$, there exists exactly one triple of $\mathcal{B}'$ containing $\{x,y\}$. There are the possible following three cases: 
\begin{enumerate}
    \item[2.1.]  $x,y \in X$;
    \item[2.2.] $x,y \in Y$;
    \item[2.3.] $x\in X, y\in Y$.
\end{enumerate}
\end{enumerate}

Case 2.1. Since $\Sigma$ is an $STS$ and $\Gamma$ does not contain pairs of vertices of $X$, there exists exactly one block in $\mathcal{B}$ containing $\{x,y\}$.

Case 2.2. There exists exactly one factor $F_j \in \mathcal{F}$ containing $x,y$. If $x_i\in X$ is such that $\varphi(x_i)=F_j$, the
triple $\{x_i, x,y\}$ is of type 2 and is it the unique, triple of $\mathcal{B}'$ containing $\{x,y\}$.

Case 2.3. Consider the factor $\varphi(x)=F_j\in \mathcal{F}$. There exists exactly one pair of $F_j$ containing $y$. Let $\{y,z\} \in F_j$. The triple $\{x,y,z\}$ is a triple of type 2 and it is the unique triple of $\mathcal{B}'$ containing $x,y$.
\end{proof}

\section{Factorization on $\mathbb Z_{2^n}$ by difference method} 

Let $r\in\mathbb N$ even and let $\mathbb Z_{2^n}=\{0,1,2,...,2^n-1\}$, $D_{2^n}=\{1,2,...,2^{n-1}\}$. In this section we define a factorization of the complete graph
$K_{2^n}$ defined in $\mathbb Z_{2^n}$, briefly a factorization on $\mathbb Z_{2^n}$, where every factor contains pair having all the same difference. Note that all the elements belong to $\mathbb Z_{2^n}$.
	
\begin{thm} \label{T:1} 
    If $a\in D_{2^n}$ is an odd number, then there exist two disjoint factors of $\mathbb Z_{2^n}$,  $F_{a,1}, F_{a,2}$, containing all the pairs $\{x,y\}\subset \mathbb Z_{2^n}$ having difference $a$.
\end{thm}
\begin{proof} 
Consider the complete graph $K_{2^n}$ defined in $\mathbb Z_{2^n}$. Since $a$ is an odd number, $pa\equiv 0 \mod 2^n$ if and only if $p\equiv 0\mod 2^n$ and so there exists in $K_{2^n}$ a unique hamiltonian cycle of length $\mathbb Z_{2^n}$, which can be described as follows:
\[ 
(0,a,2a,3a,....,2^{n}-a,0),
\]
and it is immediate to see that it can be decomposable in the following two factors:
\[
F_{a,1}:\hspace{1.5pt}\{0,a\},\{2a,3a\},....,\{2^{n}-2a,2^{n}-a\}
\]
\[
F_{a,2}:\hspace{1.5pt}\{a,2a\},\{3a,4a\},\{5a,6a\},....,\{2^{n}-a,2^{n}=0\}.
\]
\end{proof}

\begin{thm} \label{T:2} 
    If $a\in D_{2^n}-\{2^{n-1}\}$ is an even number and $k=GCD(a,2^n)$, then there exist two disjoint factors of $\mathbb Z_{2^n}$,
 $G_{a,1}, G_{a,2}$,  containing all the pairs $\{x,y\}\subset \mathbb Z_{2^n}$ having difference $a$.
\end{thm}
\begin{proof}
    In this case, since $pa\equiv 0\mod 2^n$ if and only if $p\equiv 0\mod \tfrac{2^n}{k}$, in the complete graph $K_{2^n}$ defined in $\mathbb Z_{2^n}$ there exist $k$ disjoint cycles of length $\frac{2^n}{k}$, described as follows:\\
\begin{align*}
    &C(a,1)=(0,a,2a,3a,\dots,2^{n}-a,0),\\
    &C(a,2)=(1,1+a,1+2a,1+3a,\dots,1+2^{n}-a,1+2^n=1),\\
    &\dots\\
    &C(a,k)=(k-1,k-1+a,k-1+2a,\dots,2^{n}+k-1-a,2^n+k-1=k-1)
\end{align*}
which can be decomposed in the following factors:
\begin{multline*}
    G_{a,1}\colon \{i,a+i\},\{2a+i,3a+i\},\{4a+i,5a+i\},\dots\\
    \{2^{n}-2a+i,2^{n}-a+i\},\quad \text{for }i=0,\dots,k-1,
\end{multline*}
\begin{multline*}
    G_{a,2}\colon\{a+i,2a+i\},\{3a+i,4a+i\},\{5a+i,6a+i\},\dots,\\
\{2^{n}-a+i,i\}\quad \text{for }i=0,\dots,k-1.
\end{multline*}
\end{proof}

\begin{thm} \label{T:3} 
    If $a=2^{n-1}$, then there exists an unique factor $H$ containing all the pairs $\{x,y\}\subset \mathbb Z_{2^n}$ having difference $a$.
\end{thm}
\begin{proof}
If $a=2^{n-1}$, it is immediate to verify that in the complete graph $K_{2^n}$, defined in $\mathbb Z_{2^n}$, there exist the following unique factor:
\[
H\colon\{0,2^{n-1}\},\{1, 2^{n-1}+1\}, \{2, 2^{n-1}+2\},\dots,\{2^{n-1}-1,2^{n}-1\}.
\]
\end{proof}

\noindent Note that the family $\mathcal{F}$, so defined: 
\begin{multline} \label{E:1}
    \mathcal{F}=\{H\}\cup \{F_{a,1}, F_{a,2}: a\in D_{2^n},\, a\text{ odd}\}\cup\\ \cup\{G_{a,1}, G_{a,2}: a\in D_{2^n},\, a\ne 2^{n-1},\, a\text{ even}\}
\end{multline} 
is a factorization defined in $\mathbb Z_{2^n}$, which we call the \textit{difference factorization} of $\mathbb Z_{2^n}$. Note that $a$ is an odd number in $F_{a,1}, F_{a,2}$ and it is an even number in $G_{a,1}, G_{a,2}$.

\section{A construction $v\longrightarrow 2v+7$}
Let $n\in \mathbb N$, $n\ge 3$, and let $v=2^n-7$. Note that $2^n-7$ is a positive integer number such that  for $n$ odd  $2^n-7 \equiv 1$ mod $6$ and for $n$ even $2^n-7 \equiv 3$ mod $6$. Note also that $2v+7=2^{n+1}-7$.

\begin{thm}   \label{T}
    Let $\Sigma=(X,\mathcal{B})$ be an $STS(v)$, where $v=2^n-7, n\in \mathbb N$, $n\ge 3$. It is possible to define an $STS(2v+7)$
$\Sigma'=(X', \mathcal{B}')$ embedding $\Sigma$.
\end{thm}
\begin{proof}
Let $\Sigma=(X,\mathcal{B})$ be an $STS(v)$ defined on $X=\{x_1, x_2, .... , x_v\}$. Further, let:
\begin{itemize}
    \item $Y=\mathbb Z_{2^n}$ be a set such that $X\cap Y=\emptyset$;
    \item $\mathcal{F}$ be the difference factorization of $\mathbb Z_{2^n}$ given in \eqref{E:1};
    \item $\{a,b,c\}$ be a difference triple contained in $D_{2^n}$, for $a\neq 2^{n-1}, b\neq 2^{n-1}, c\neq 2^{n-1}$, with either $c=a+b$ or $a+b+c=2^n$;
    \item $\mathcal{F}^*$ be the family of factors obtained by $\mathcal{F}$ by excluding the factors containing the pairs having differences either $a$ or $b$ or $c$.
\end{itemize}
Further, since $|\mathcal{F}^*|=|\mathcal{F}|- 6 = 2^n -7 = v$, it is possible to define any bijection $\varphi: X\longrightarrow \mathcal{F}^*$ from $X$ into $\mathcal{F}^*$.  

In what follows we fix $a=1, b=2, c=3$. Note that, from Theorems \ref{T:1}, \ref{T:2}, \ref{T:3}, this does not harm the generality. At this point, define the hypergraph $\Sigma'=(X',\mathcal{B}')$ as follows:
\begin{itemize}
    \item $X'=X\cup Y$;
\item $\mathcal{B}'=\mathcal{B}\cup \Gamma \cup \Delta$, where:
    \[
\Gamma=\{ \{x,y,z\}: x\in X, \{y,z\}\in \varphi(x)\},
\]
\[
\Delta=\{ \{i,i+1,i+3\}: i=0,1,2,...,2^n-1\}.
\]
\end{itemize}
We call the triples of $\mathcal{B}$ $\textit {of type 1}$, the triples of $\Gamma$ \textit {of type 2},
the triples of $\Delta$ \textit {of type 3}. Observe that the triples of $\Delta$ are all the translates generated by the base block $\{0,1,3\}$, defined by the difference triple $\{1,2,3\}$. We are going to prove that $\Sigma'$ is an $STS(2v+7)$. Indeed:
\begin{enumerate}
    \item it is immediate that $\Sigma'$ is a hypergraph of order $|X'|=2v+7$, uniform of rank $3$;
    \item for every $x,y\in X', x\neq y$, there exists exactly one triple of $\mathcal{B}'$ containing ${x,y}$. Consider the following three possible cases:
	\begin{enumerate}
	    \item[2.1.] $x,y \in X$;
	    \item[2.2.1.] $x,y \in Y$ and the difference between $x,y$ is either $1$ or $2$ or $3$;\\
	    \item[2.2.2.] $x,y \in Y$ and the difference between $x,y$ is in $\{4,\dots,2^{n-1}\}$;\\
	    \item[2.3.] $x\in X, y\in Y$.
	\end{enumerate}
\end{enumerate}

Case 2.1. Since $\Sigma$ is an $STS$ and $\Gamma \cup \Delta$ do not contain vertices of $X$, there exists exactly one block in $\mathcal{B}$, of type $1$, containing $x,y$.

Case 2.2.1. In $\Sigma$ there are not vertices of $Y$ and in the factors of $\mathcal F^*$ there are not pairs of vertices having differences $1$, $2$, or $3$. In $\Delta$ there are triples containing all the pairs of vertices of $Y$ with differences either $1$ or $2$ or $3$, therefore there exists in $\Delta$ a block, of type $3$, containing $\{x,y\}$.

Case 2.2.2. In the factors of $\mathcal F^*$ there are all the pairs of $Y$ having difference $a\in D_{2^n} - \{1,2,3\}$, therefore there exists in $\mathcal F^*$ a factor $F$ containing $\{x,y\}$. If $z=\varphi^{-1}(F)$, then the triple $\{x,y,z\}$ is a bock of type $2$ of $\Sigma'$.

Case 2.3. In this case, consider the factor $\varphi(x)=F\in \mathcal F^*$. There exists exactly one pair of the factor $F$ containing $y$. Let $\{y,z\} \in F$. The triple $\{x,y,z\}$ is a triple of type 2 containing the pair $\{x,y\}$.

So, we have proved that, in general, for every pair $\{x,y\}\subset X'=X\cup Y$ it is $d(x,y)\geq 1$. To prove that it is exactly $d(x,y)=1$ it is sufficient to prove that $\mathcal{B}'=\frac{(2v+7)(2v+6)}{6}$, which is the exact number of blocks in any 
$STS(2v+7)$. Since $\mathcal{B}'=\mathcal{B}\cup \Gamma \cup \Delta$ and $\mathcal{B}, \Gamma, \Delta$ are pairwise disjoint, it follows that:
\[
|\mathcal{B}'|=|\mathcal{B}|+|\Gamma|+|\Delta|=\frac{v(v-1)}{6}+v \cdot \frac{v+7}{2}+(v+7)=\frac{(2v+7)(2v+6)}{6}.
\]
\end{proof}

\begin{coro}  \label{C}
    Let $n\in \mathbb N$, $n\ge 4$, and let $\Sigma=(X,\mathcal B)$ be an $STS(v)$, with $v=2^n-7$, having a maximal independent set $T$ of maximal cardinality $\tfrac{v-1}{2}$. Then there exists an $STS$ $\Sigma'$ of order $v'=2^{n+1}-7$ embedding $\Sigma$ with a maximal independent set $T'\supset T$ of maximal cardinality $\tfrac{v'-1}{2}$.
\end{coro}
\begin{proof}
    In the proof of Theorem \ref{T} if $T$ is a maximal independent set for $\Sigma$, then $T\cup\{2i\mid i=0,1,\dots,2^{n-1}\}$ is a maximal independent set for $\Sigma'$ of cardinality $\tfrac{v'-1}{2}$, with $v'=2^{n+1}-7$, under the following conditions:
\begin{itemize}

    \item the restriction $\varphi|_T\colon T\rightarrow \mathcal F^{\star}$ induces a bijection with the subset of $\mathcal F^*$ determined by the differences $a\in D_{2^n}$, $a$ odd;
    \item the restriction $\varphi|_{X-T}\colon X-T\rightarrow \mathcal F^{\star}$ induces a bijection with the subset of $\mathcal F^*$ determined by the differences $a\in D_{2^n}$, $a$ even.

\end{itemize} 
\end{proof}

\begin{coro}
    There exists for any $n\in\mathbb N$, $n\ge 4$, an $STS(v)$, with $v=2^n-7$, having a maximal independent set of maximal cardinality $\tfrac{v-1}{2}$.  
\end{coro}
\begin{proof}
It is sufficient to apply iteratively the previous corollary, considering that the base case is $v=9$, for which the easy statement is proved in the Appendix.  
\end{proof}

\begin{coro}
    Let $n\in \mathbb N$, $n\ge 4$. Then there exists an STS$(v)$ $\Sigma$, with $v=2^n-7$, such that $\overline{\chi}(\Sigma)=n-1$.
\end{coro}
\begin{proof}
    Note that by \cite[Theorem 1 and Corollary 1]{MT} for any bicolorable $STS(2^n-7)$ $\Sigma$ it must be $\overline{\chi}(\Sigma)<n$, which implies that we simply need to prove the existence of an $STS(2^n-7)$ which is $(n-1)$-bicolorable.

    Let $v=9$. In this case, given the system $\Sigma$ in the Appendix, it is $3$-bicolorable with the color classes $\{x_1,x_2,x_4,x_5\}$, $\{x_3,x_6,x_7,x_8\}$ and $\{x_9\}$. Now, let $T=\{x_1,x_2,x_4,x_5\}$ a maximal independent set for $\Sigma$. When we apply the construction given in Corollary \ref{C} to $\Sigma$, the system $\Sigma'$ of order $25$ that we get is $4$-bicolorable, since we can take as color classes $C_1=\{x_1,x_2,x_4,x_5\}\cup \{2i\mid i=0,\dots,7\}$, $C_2=\{x_3,x_6,x_7,x_8\}$, $C_3=\{x_9\}$ and $C_4=\{2i+1\mid i=0,,\dots,7\}$. Moreover, a maximal independent is $T\cup \{2i\mid i=0,\dots,7\}=C_1.$

    For $n\ge 5$ we apply iteratively the construction given in Corollary \ref{C}, where the case $n=5$ has been previously explained. So, we have an $STS(2^{n-1}-7)$ which is $(n-2)$-bicolorable and such that there exists a $(n-2)$-bicoloring having as color class, say $C_1$, a maximal independent set $T$. In this case, proceeding as in Corollary \ref{C} it is sufficient to give the color of $1$ to the vertices in $\{2i\mid i=0,\dots,2^{n-1}-1\}$ and the color $n-1$ to the vertices in $\{2i+1\mid i=0,\dots,2^{n-1}-1\}$. Note that in this way a color class in an $(n-1)$-bicoloring of $\Sigma'$ coincides with a maximal independent set of maximal cardinality. 
\end{proof}

\section{Appendix}
In this section we give an application of the construction described in this paper in the case $v=9$ and $2v+7=25$. Let:
\begin{itemize}
    \item $Y=\mathbb Z_{16}=\{0,1,2,3,4,5,6,7,8,9,10,11,12,13,14,15\}$,
    \item $D_{16}=\{1,2,3,4,5,6,7,8\}$,
    \item $\Sigma=(X,\mathcal{B}), STS(9)$,
    \item $X=\{x_1,x_2,x_3,x_4,x_5,x_6,x_7,x_8,x_9\}$,
    \item $\mathcal{B}:$
	\[
	\{x_1,x_2,x_3\},\hspace {2.00pt}\{x_1,x_4,x_7\},\hspace {2.00pt}\{x_1,x_5,x_9\},\hspace {2.00pt}\{x_1,x_6,x_8\},
    \]
    \[    
\{x_4,x_5,x_6\},\hspace {2.00pt}\{x_2,x_5,x_8\},\hspace {2.00pt}\{x_2,x_6,x_7\},\hspace {2.00pt}\{x_2,x_4,x_9\},
\]
\[
\{x_7,x_8,x_9\},\hspace {2.00pt}\{x_3,x_6,x_9\},\hspace {2.00pt}\{x_3,x_4,x_8\},\hspace {2.00pt}\{x_3,x_5,x_7\}.
\]
\end{itemize}
Note that $T=\{x_1,x_2,x_4,x_5\}$ is a maximal independent set for $\Sigma$ and that the sets $\{x_1,x_2,x_4,x_5\}$, $\{x_3,x_6,x_7,x_8\}$ and $\{x_9\}$ are the color classes of a $3$-bicoloring of $\Sigma$. Let $\Gamma$ be the family of triples containing an $x_i \in X$ and a pair indicated in the follow columns:

\begin{minipage}[t]{0.16\textwidth}
$x_1$\\
$\downarrow$\\
$F_{5,1}:$\\
$0\cdot 5$\\
$10 \cdot 15$\\
$4 \cdot 9$\\
$14 \cdot 3$\\
$8 \cdot 13 $\\
$2 \cdot 7$\\
$12 \cdot 1$\\
$6 \cdot 11$\\
\end{minipage}
\begin{minipage}[t]{0.16\textwidth}
$x_2$\\
$\downarrow$\\
$F_{5,2}:$\\
$5\cdot 10$\\
$15 \cdot 4$\\
$9 \cdot 14$\\
$3 \cdot 8$\\
$13 \cdot 2 $\\
$7 \cdot 12$\\
$1 \cdot 6$\\
$11 \cdot 0$\\
\noindent
\end{minipage}
\begin{minipage}[t]{0.16\textwidth}
$x_3$\\
$\downarrow$\\
$G_{4,1}:$\\
$0\cdot 4$\\
$8 \cdot 12$\\
$1 \cdot 5$\\
$9 \cdot 13$\\
$2 \cdot 6 $\\
$10 \cdot 14$\\
$3 \cdot 7$\\
$11 \cdot 15$\\
\end{minipage}
\begin{minipage}[t]{0.16\textwidth}
$x_4$\\
$\downarrow$\\
$F_{7,2}:$\\
$7\cdot 14$\\
$5 \cdot 12$\\
$3 \cdot 10$\\
$1 \cdot 8$\\
$15 \cdot 6 $\\
$13\cdot 4$\\
$11 \cdot 2$\\
$9 \cdot 0$\\
\end{minipage}

\begin{minipage}[t]{0.16\textwidth}
$x_5$\\
$\downarrow$\\
$F_{7,1}:$\\
$0\cdot 7$\\
$14 \cdot 5$\\
$12 \cdot 3$\\
$10 \cdot 1$\\
$8 \cdot 15 $\\
$6 \cdot 13$\\
$4 \cdot 11$\\
$2 \cdot 9$\\
\end{minipage}
\begin{minipage}[t]{0.16\textwidth}
$x_6$\\
$\downarrow$\\
$G_{4,2}:$\\
$4\cdot 8$\\
$12 \cdot 0$\\
$5 \cdot 9$\\
$13 \cdot 1$\\
$6 \cdot 10$\\
$14\cdot 2$\\
$7 \cdot 11$\\
$15 \cdot 3$\\
\end{minipage}
\begin{minipage}[t]{0.16\textwidth}
$x_7$\\
$\downarrow$\\
$F_{6,1}:$\\
$0\cdot 6$\\
$12 \cdot 2$\\
$8 \cdot 14$\\
$4 \cdot 10$\\
$1 \cdot 7$\\
$13\cdot 3$\\
$9 \cdot 15$\\
$5 \cdot 11$\\
\end{minipage}
\begin{minipage}[t]{0.16\textwidth}
$x_8$\\
$\downarrow$\\
$F_{6,2}:$\\
$6\cdot 12$\\
$2 \cdot 8$\\
$14 \cdot 4$\\
$10 \cdot 0$\\
$7 \cdot 13 $\\
$3 \cdot 9$\\
$15 \cdot 5$\\
$11 \cdot 1$\\
\end{minipage}
\begin{minipage}[t]{0.16\textwidth}
$x_9$\\
$\downarrow$\\
$H:$\\
$0 \cdot 8$\\
$1 \cdot 9$\\
$2 \cdot 10$\\
$3 \cdot 11$\\
$4 \cdot 12$\\
$5 \cdot 13$\\
$6 \cdot 14$\\
$7 \cdot 15$\\
\end{minipage}

\noindent Let $X'=X\cup Y$ and $\mathcal{B}'=\mathcal{B}\cup \Gamma \cup \Delta$, where $\mathcal{B}$ and $\Gamma$ are already defined and  
\[
\Delta=\{ \{i,i+1,i+3\}: i=0,1,2,\dots,15\}.
\]
It is immediate to verify that $\Sigma'(X',\mathcal{B}')$ is an $STS(25)$.

\vspace{1cc}

\end{document}